\documentclass[twocolumn]{svjour3}

\usepackage{amsmath}
\usepackage{amssymb}
\usepackage{graphics}
\usepackage{epsfig}
\usepackage{times}
\usepackage{prettyref}
\usepackage[numbers,sort]{natbib}
\title{
Using ODE waveform-relaxation methods to efficiently include gap junctions in distributed neural network simulations
}

 \author{Matthias Bolten \and Jan Hahne}
 \institute{Matthias Bolten \and Jan Hahne \at School of Mathematics and Natural Sciences \\University of Wuppertal \\ Wuppertal, Germany \\\email{hahne@math.uni-wuppertal.de}}

\newrefformat{fig}{Fig.~\ref{#1}}
\newrefformat{theo}{Theorem~\ref{#1}}

\begin{document}

\maketitle
\thispagestyle{empty}
\pagestyle{empty}

%%%%%%%%%%%%%%%%%%%%%%%%%%%%%%%%%%%%%%%%%%%%%%%%%%%%%%%%%%%%%%%%%%%%%%%%%%%%%%%%
\begin{abstract}

Waveform-relaxation methods divide systems of differential equations into subsystems and therefore allow for parallelization across the system. Here we present an application for ODE waveform-relaxation methods in the context of spiking neural network simulators. Parallel spiking neural network simulators make use of the fact that the dynamics of neurons with chemical synapses is decoupled for the duration of the minimal network delay and thus can be solved independently for this duration. The inclusion of electrical synapses, so-called gap junctions, requires continuous interaction between neurons and therefore constitutes a conceptional problem for those simulators. We present a suitable waveform-relaxation method for an efficient integration of gap junctions and demonstrate that the use of the waveform-relaxation method improves both, accuracy and performance, compared to a non-iterative solution of the problem. We investigate the employed method in a reference implementation in the parallel spiking neural network simulator NEST.

\end{abstract}

%%%%%%%%%%%%%%%%%%%%%%%%%%%%%%%%%%%%%%%%%%%%%%%%%%%%%%%%%%%%%%%%%%%%%%%%%%%%%%%%
\section{Introduction}

\begin{sloppypar}
Spiking neural network simulations play an important role in neuroscience as they allow researchers to study the behavior and properties of large neural networks and  to investigate new models and hypotheses on a biologically realistic scale.
In order to provide this capability for researchers the corresponding simulators need to be optimized for parallel computation on supercomputers. 
This highly optimized structure, however, restrains the efficient inclusion of new features, if those features violate the assumptions the optimized structure is based on. One example for such a case is the inclusion
of electrical synapses, so-called \emph{gap junctions}. The computational and neuroscientific aspects of the 
implementation of gap junctions in a spiking neural network simulator (without any special focus on how waveform-relaxation methods can be used to overcome the difficulties of including gap junctions) are described in \citep{Hahne15_00022, Hahne16}. Here we discuss the missing aspects with regard to the waveform-relaxation method and show how exactly they can improve the inclusion of gap junctions.
After a brief introduction to spiking neural network simulators and waveform-relaxation methods for ODEs, \prettyref{sec:including_gj} demonstrates that the use of the waveform-relaxation method improves both, accuracy and performance, 
compared to a non-iterative way to include gap junctions. Finally \prettyref{sec:conclusion} gives some concluding remarks.
\end{sloppypar}

\subsection{Spiking neural network simulators}
One strategy to describe networks of neurons in computational neuroscience is by modeling them using a spiking neural network model.
This approach is motivated by the microscopic dynamics of individual neurons and is therefore called a \emph{bottom-up approach}.
The idea is to model the individual neuron dynamics and their interaction with each other. The individual neurons can be described as
\emph{hybrid systems}:
\begin{eqnarray}
\frac{dy}{dt} & = &f(y(t)) \label{eq:hybrid1} \\
y(t) & \leftarrow  & g_j(y(t))  \quad \text{upon spike from synapse } j \label{eq:hybrid2}
\end{eqnarray}

\begin{sloppypar}
\noindent Here the state $y(t)=(y_{1}(t),\ldots,y_{i}(t),\ldots,y_{N}(t))^T$ of the neuron evolves continuously to some biophysical equations \prettyref{eq:hybrid1} described by the $N$-dimensional function $f:\left[t_0 , t_0+T \right] \times \mathbb{R}^{N}\rightarrow\mathbb{R}^{N}$ and a spike
received by the neuron triggers a instantaneous change in some state variable $y_i$ \prettyref{eq:hybrid2} described by synapse-specific functions $g_j:\left[t_0 , t_0+T \right] \times \mathbb{R}^{N}\rightarrow\mathbb{R}^{N}$. The neuron emits a
spike if its membrane potential $V$ satisfies some threshold condition, e.g.\ $V > \theta$.
Typically the membrane potential is taken as the first state variable, i.e. $y_1=V$.
\end{sloppypar} 
 
The focus of spiking neural network simulators (for a review see \citep{Brette07_349}) ranges from detailed
neuron morphology (\textit{NEURON} \citep{Carnevale06}, \textit{GENESIS} \citep{Bower07_1383})
to an abstraction of neurons without spatial extent (\textit{NEST} \citep{Gewaltig_07_11204},
\textit{BRIAN} \citep{Goodman13}). Here we will focus on the latter. 
In this case, the biophysical equations \prettyref{eq:hybrid1} of the individual neurons are ODEs, only dependent on the time $t$.
Due to the missing spatial dimensions, the neurons in this modeling approach are also called \emph{point neurons}.

In networks of point neurons each chemical synapse has its own synaptic delay $d_{ij}$, which depends on the emitting neuron $i$ and the receiving neuron $j$. 
These delays are the result of different biological processes which cannot be modeled more accurately for point neurons, e.g.\ the time the
postsynaptic potential needs to travel from the synapse
on a dendrite to the soma. Thus, a spike which is emitted at time $t_1$ in neuron $i$ affects the connected neuron $j$
at time $t_1 + d_{ij}$. In a network, the dynamics of all neurons is decoupled for the duration of the minimal
network delay $d_{\mathrm{min}}=\min_{i,j}(d_{ij})$.  Hence, the dynamics of each neuron can be propagated
independently for the duration $d_{\mathrm{min}}$ without requiring
information from other neurons.

\begin{figure}[!tbh]
\begin{centering}
\includegraphics[width=0.47\textwidth]{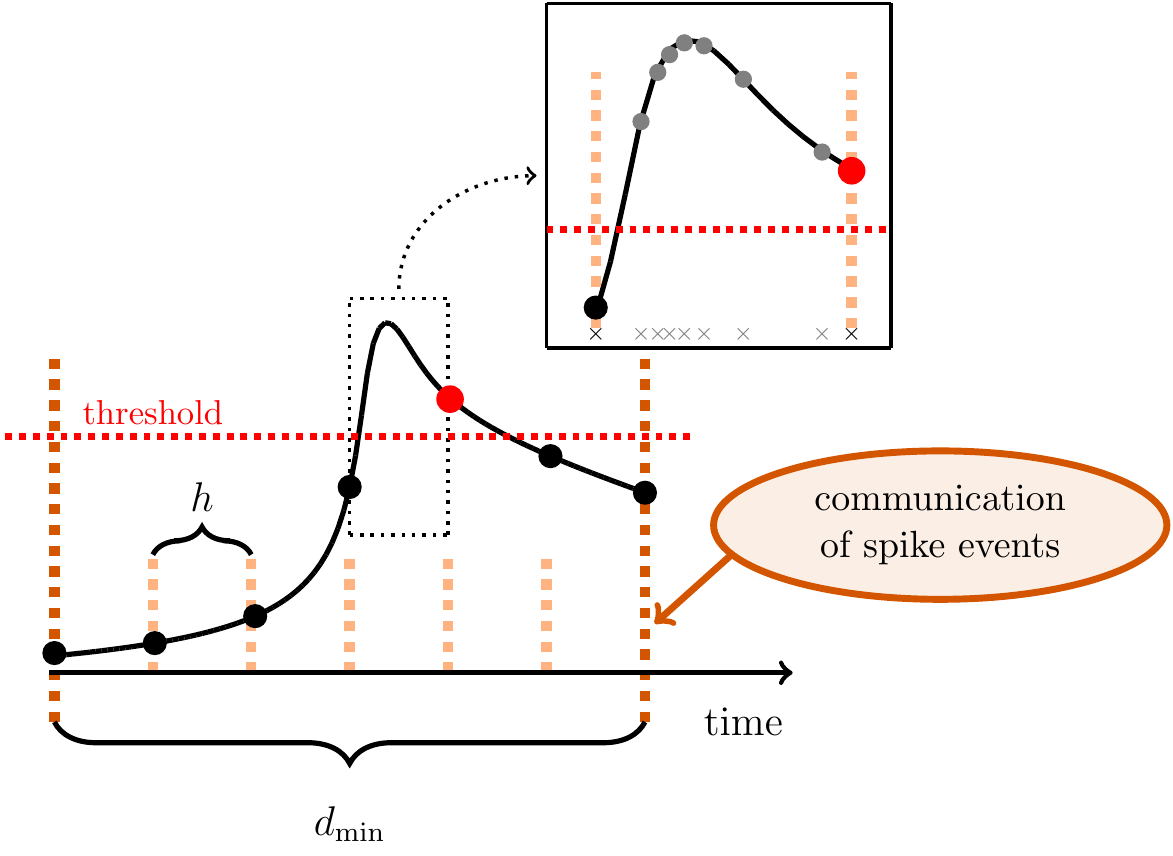}
\par\end{centering}
\caption[Time course of the membrane potential of a single neuron during one $d_{\min}$-interval]{\textbf{Time course of the membrane potential of a single neuron during one $d_{\min}$-interval.} Black dots mark the values of the membrane potential at the grid points. The red dot indicates the point in time where the spike triggered by the threshold crossing is registered.
The enlarged rectangle shows the gray marked internal grid points used by the model specific numerical method solving the neuron dynamics. 
\label{fig:setup-solver}}
\end{figure}

Efficient simulators make use of the delayed and point-event like nature of the spike interaction by distributing neurons across available 
processes and communicating spikes only after this period \citep{Morrison05a}. \prettyref{fig:setup-solver} shows the progress of the widely used NEST simulator over one interval of duration of the minimal network delay,
from the single neuron perspective. All times are restricted to a time grid with fixed step size $h$, which means that every spike time and all information on the current state of a neuron is restricted to this grid, 
and that every delay of a synaptic connection has to be a multiple of $h$.

\vspace*{-0.4cm}

\noindent 

A numerical method solving the single neuron dynamics \prettyref{eq:hybrid1} is thus required to produce approximations of the state 
of the neuron at each grid point. This does, however, still allow the use of a method with adaptive step size within one interval of duration $h$. 
The state of the neuron at those additional evaluation points within the interval cannot, however, be recorded by any NEST device.

\subsection{Waveform-relaxation methods}

Waveform-relaxation methods are a set of iterative methods for the solution of
systems of differential equations. The name and the concept were first introduced in the early 1980s for the simulation of large-scale electric circuits \citep{Lelarasmee82,Lelarasmee82_1}. 
For any given $N$-dimensional system of ODEs $y^{\prime}(t)=f(t,y(t))$ with $f:\left[t_0 , t_0+T \right] \times \mathbb{R}^{N}\rightarrow\mathbb{R}^{N}$ and initial value $y(t_{0})=y_0$,
the basic idea is to divide the ODE-system into $v \leq N$ preferably weakly
coupled subsystems

\begin{equation}
y_{i}^{\prime}(t)  =  f_{i}(t,y_{1}(t),\ldots,y_{v}(t)) \quad \; i=1,\ldots,v
\end{equation}

\noindent  and to solve each subsystem independently by treating the influence
of the other subsystems as given input. For $v < N$ the resulting method is called a \emph{block method}. 

Starting with an initial guess $y_{i}^{(0)}(t)$ for the solution of
each subsystem over the entire iteration interval $[t_{0},t_{0}+\mathcal{T}]$,
the solution of the original ODE-system is determined by iteratively
solving the independent subsystems, where the input for the $i$-th subsystem $y_{1}(t),...,y_{i-1}(t),y_{i+1}(t),...,y_{v}(t)$ is based on previously obtained
solutions from the current or any of the previous iterations and hence acts as a
given input to the $i$-th system. In the absence of other reasonable extrapolations the initial guess
$y_{i}^{(0)}(t)$ is usually chosen to be constant as $y_{i}^{(0)}(t)=y_{0}$ for all $t\in[t_{0},t_{0}+\mathcal{T}]$ in order to fulfill the initial value condition.
Two prominent examples of waveform-relaxation methods are 

\begin{itemize}
\item[i)] the (block) Jacobi waveform-relaxation method

\begin{equation}
\begin{split}
{y_{i}^{\prime}}^{(m)}(t) =  f_{i}(t,y_{1}^{(m-1)}&(t),...,y_{i-1}^{(m-1)}(t),y_{i}^{(m)}(t), \\ & y_{i+1}^{(m-1)}(t),...,  y_{v}^{(m-1)}(t)) \\
\end{split}
\label{eq:wfr-jacobi-version}
\end{equation}

with $i=1,\ldots,v$, where the input for the $m$-th iteration is
completely based on the solutions of the $(m-1)$-th iteration. For each iteration
this strategy enables parallel processing of all subsystems and is
hence well-suited for distributed simulations.

\item[ii)] the (block) Gauss-Seidel waveform-relaxation method

\begin{equation}
\begin{split}
{y_{i}^{\prime}}^{(m)}(t) =  f_{i}(t,y_{1}^{(m)}&(t),...,y_{i-1}^{(m)}(t), y_{i}^{(m)}(t),\\ &y_{i+1}^{(m-1)}(t),...,y_{v}^{(m-1)}(t)) \\
\end{split}
\end{equation}

with $i=1,\ldots,v$, where the input of subsystem $i$ is taken from the current iteration for all systems $j<i$ and from the previous iteration
for the remaining subsystems. Compared to the Jacobi version, this strategy is expected to speed up the convergence, but comes at the price
of less potential for parallel processing.
\end{itemize}

\begin{sloppypar}
Waveform-relaxation methods can be seen as an extension of the Picard-Lindel\"{o}f iteration
which was first mentioned almost one century before the first waveform-relaxation methods were developed \citep{Lindeloef1894}. 
Therefore some works, e.g.\ \citep{Nevanlinna89_1,Nevanlinna89_2}, also refer to waveform-relaxation methods as the (generalized) \emph{Picard-Lindel\"{o}f iteration}.
Another quite commonly used term is \emph{dynamic iteration}, e.g.\ in \citep{Bjorhus94,Bjorhus95}. This name explores the common ground between waveform-relaxation methods for
linear ODEs and iterative methods for the solution of linear systems of algebraic equations (e.g.\ the Jacobi or Gauss-Seidel method), 
which in this view can be seen as the \emph{static iteration} counterpart. 
Due to their setup, waveform-relaxation methods are particularly interesting for parallel simulations of large-scale problems. 
In the class of parallelizable methods, they can be classified as \emph{parallel across the system} (see \citep{Burrage93} for a review of parallelizable methods for ODEs).
\end{sloppypar}

Waveform-relaxation methods for ODEs and their convergence have been studied by various authors (for a quite extensive list see the introduction of \citep{Fan11}).
In general one can distinguish between studies of \emph{continuous} and \emph{discrete} waveform-relaxation methods. The former investigate the 
application of the waveform-relaxation method directly to the continuous problem while the latter first apply some numerical method to the continuous formulation of the problem and investigate the convergence and behavior of the resulting waveform-relaxation method
in comparison with the solution of the applied numerical scheme without the use of the waveform-relaxation technique.

\begin{sloppypar}
In both cases, a formulation of the method with a so-called \emph{splitting function} $F:\left[t_0 , t_0+T \right] \times \mathbb{R}^{N} \times \mathbb{R}^{N} \rightarrow\mathbb{R}^{N}$ with
\begin{equation}
F(t,y,y)=f(t,y) \quad  \forall t \in \left[t_0 , t_0+T \right] \text{ and } \forall y \in \mathbb{R}^{N} \,.
\end{equation}
allows the simultaneous analysis of different waveform-relaxation methods. One thus considers the problem
\begin{equation}
{y^{\prime}}^{(m)}(t)=F(t,y^{(m)}(t),y^{(m-1)}(t))\,.
\label{eq:wfr-with-splitting}
\end{equation}
\end{sloppypar}

\begin{sloppypar}
The choice $F(t,y,z)=f(t,z)$ yields the Picard-Lindel\"{o}f iteration, 
while $F_i(t,y,z) = f_{i}(t,z_{1},...,z_{i-1},y_{i},z_{i+1},...,z_{v})$ results in the block Jacobi ($v<N$) or Jacobi ($v=N$) waveform-relaxation method
\prettyref{eq:wfr-jacobi-version}. 
% The universal formulation also includes non-standard waveform-relaxation methods where a particular part of a function $f$ containing $y_i(t)$ is evaluated 
% with information from the current iteration, and the remaining part is evaluated with information from the previous iteration. For example, for the one-dimensional problem $f(t,y)=y^2$ the 
% splitting function $F(t,y,z) = 0.7 y^2 + 0.3 z^2$ defines such a non-standard method.
The superlinear convergence of continuous waveform-relaxation methods can be shown with the following theorem that has been stated by many authors (e.g.\ in \citep{Bellen93,Bjorhus95,Burrage93,intHout95,Nevanlinna89_1}).
\end{sloppypar}

\begin{theorem}[Superlinear convergence]
Let the splitting function $F$ satisfy
\begin{equation*}
\begin{split}
 \lVert F(t,x,y)-F(t,\hat{x},y) \rVert_2 & \,\; \leq  \, K_1 \lVert x-\hat{x} \rVert_2 \\ & \forall x,\hat{x},y \in \mathbb{R}^{N},\; \forall t \in \left[t_0 , t_0+\mathcal{T} \right]
\end{split}
\end{equation*}
and
\begin{equation*}
\begin{split}
 \lVert F(t,x,y)-F(t,x,\hat{y}) \rVert_2 & \,\; \leq  \, K_2 \lVert y-\hat{y} \rVert_2 \\ & \forall x,y,\hat{y} \in \mathbb{R}^{N},\; \forall t \in \left[t_0 , t_0+\mathcal{T} \right]
\end{split}
\end{equation*}

\begin{sloppypar}\noindent for some constants $K_1$ and $K_2$. Then the sequence $y^{(m)}(t), \; m=0,1,\ldots$ produced by the waveform-relaxation method \prettyref{eq:wfr-with-splitting} converges superlinearly to
the solution $y(t)$ of the corresponding initial value problem, such that\end{sloppypar}
\begin{equation}
\begin{split}
& \sup\limits_{t_0 \leq t \leq t_0 + \mathcal{T}}   \lVert y^{(m)}(t) - y(t) \rVert_2 \\ & \quad \quad \leq e^{K_1\mathcal{T}} 
\frac{(K_2\mathcal{T})^{m}}{m!} \sup\limits_{t_0 \leq t \leq t_0 + \mathcal{T}} \lVert y^{(0)}(t) - y(t) \rVert_2 \,.
\label{eq:cont-wfr-conv-rate}
\end{split}
\end{equation}
\label{theo:bjorhus-continuous}
\end{theorem}

\section{Including gap junctions in a spiking neural network simulator}\label{sec:including_gj}

\begin{figure}[!tbh]
\begin{centering}
\includegraphics[width=0.47\textwidth]{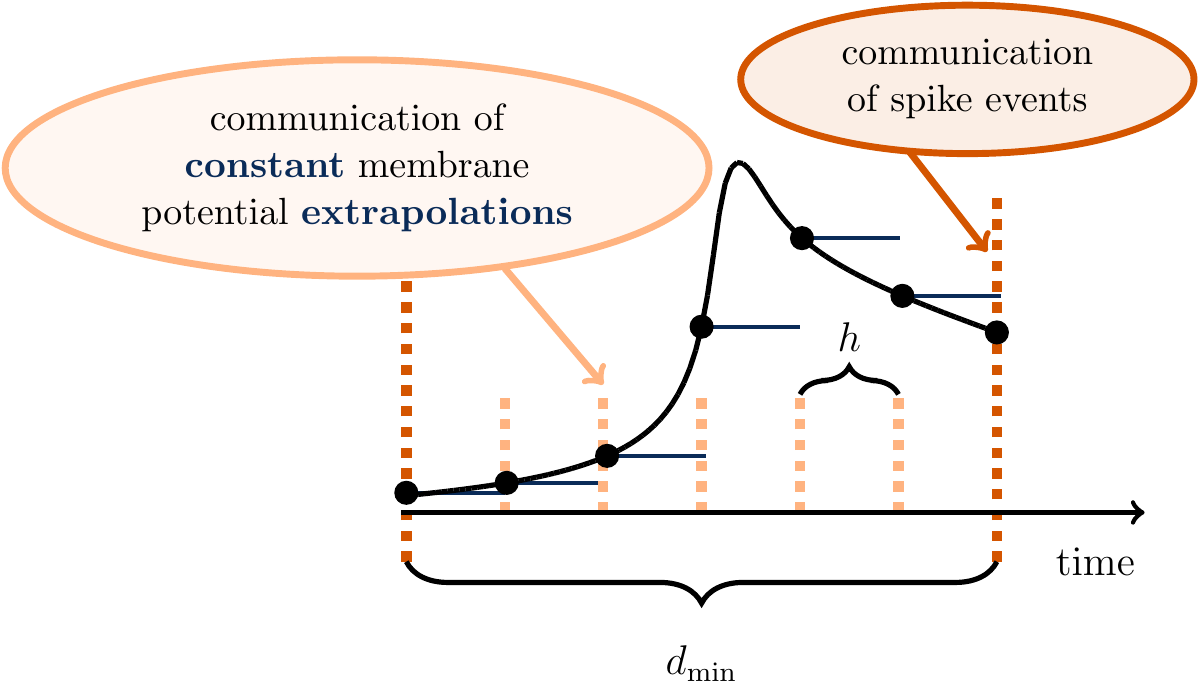}
\par\end{centering}
\caption[Progress of the non-iterative solution during one $d_{\min}$-interval]{\textbf{Progress of the non-iterative solution during one $d_{\min}$-interval.} 
Black dots mark the values of the membrane potential at the grid points, which are communicated to the connected neurons after every interval of duration $h$.  
\label{fig:straight-forward}}
\end{figure}

\begin{sloppypar}
Recently, advances in molecular biology
revealed the widespread existence of gap junctions in the mammalian nervous system,
which suggests their diverse roles in learning and memory, movement
control, and emotional responses \citep{Connors04,Hormuzdi04_113,Dere11_206}. This makes them a desirable biological feature for spiking neural network simulators. 
However, their efficient inclusion constitutes a non-trivial problem, as gap junctions require instantaneous interaction between neurons.
Gap junctions are typically represented by a current of the form
\begin{equation}
\begin{split}
I_{\mathrm{gap},i}(t)  =  g_{ij}\,(V_{j}(t)&-V_{i}(t)) \\ & \text{and} \quad I_{\mathrm{gap},j}(t)=-I_{\mathrm{gap},i}(t)  \label{eq:Igap} \;.
\end{split}
\end{equation}
occuring in both cells at the site of the gap
junction. In point-neuron models this current immediately affects the membrane potential.
As a result of a gap-junction connection between neurons $i$ and $j$ the originally decoupled hybrid systems of these
neurons are combined as an interdependent system of ODEs. Any additional gap-junction 
connection between another neuron and either $i$ or $j$ adds a further set of equations
to the coupled system. In a biologically realistic simulation of a
local cortical network, each neuron has dozens of gap-junction
connections. Consequently, the dynamics of almost all neurons are
likely interrelated by one large system of ODEs. Here the optimized structure of the simulators for spiking synapses constitutes a problem, as any communication between neurons involves communication across processes and possibly across compute nodes and is only possible on the grid with step size $h$.
Without loss of generality we employ an explicit Runge-Kutta-Fehlberg method with adaptive step-size control for the solution of the single-neuron dynamics \prettyref{eq:hybrid1} within one time step of length $h$ and use a Hodgkin-Huxley neuron model for our tests (for a general mathematical introduction see \citep{Boergers17}, for the model equations see \citep[Section 5.2]{Hahne18}).
\prettyref{fig:straight-forward} shows the only reasonable non-iterative solution to the problem over one interval of the duration of the minimal network delay $d_{\min}$ from the single neuron perspective: 
the membrane potentials of gap-junction coupled neurons are
communicated at the beginning of each $h$-interval and are assumed to be constant for the duration of that
time step. Note that the solver may still use several steps to cover the interval of duration $h$.
\end{sloppypar}

In general this constant extrapolation is the best guess, as the future evolution of the membrane potential is unknown and 
can be affected either by excitatory or by inhibitory inputs to the neuron.
However, in the time during a spike the future evolution of the membrane potential can be estimated in a more precise way: 
In a Hodgkin-Huxley point-neuron model 
a spike is triggered by a sort of chain reaction of the channel dynamics that occurs when the membrane potential increases. 
This chain reaction forms the shape of an action-potential in the membrane potential. Given a fixed set of model parameters, the temporal evolution of the membrane potential during and shortly after a spike is the same for every occurrence of a spike and can therefore be predicted in advance. During this period of time the membrane potential reaches values above a parameter-determined threshold that cannot be reached without the occurrence of a spike. Therefore a membrane potential above this threshold is a safe indicator for the occurrence of a spike in the near future or past.

   \begin{figure}[!tbh]
      \centering
      \includegraphics{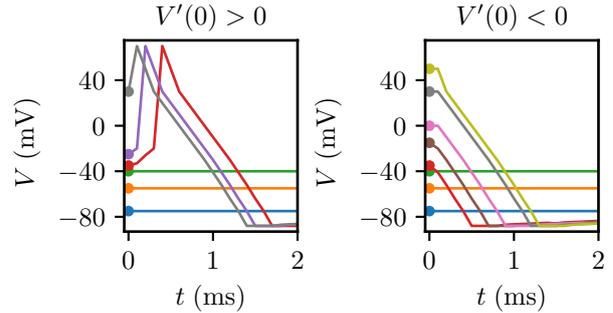}
      \caption[Improved extrapolation strategy with spike detection]{\textbf{Improved extrapolation strategy with spike detection.} The colored curves show the chosen extrapolation for the accordingly colored initial values of the membrane potential $V$ marked with circles. The figure assumes a parameter-determined spike threshold of $-40\,\mathrm{mV}$.}
      \label{fig:spike-detection}
   \end{figure}

   \begin{sloppypar}
\prettyref{fig:spike-detection} shows how this information can be used to improve the extrapolation strategy. 
Whenever the membrane potential is above the threshold the derivative $V^{\prime}(t_0) = f_1(t_0)$ is evaluated to investigate if the membrane potential is
still on the rise or is already decreasing again. Then the extrapolation is chosen based on a precomputed approximation of the uniform shape of an action-potential: As both increase and decrease are strictly monotone we can easily find the correct starting point for our approximation by evaluating $V(t_0)$ and mapping it to the precomputed approximation. If the membrane potential is below the threshold we still use a constant extrapolation.
We call this additional extrapolation strategy \emph{spike detection}.
\end{sloppypar}

We investigate the non-iterative method with a simple two-neuron network: two identical model
neurons $i$ and $j$ that are coupled by a gap junction should behave exactly the
same as an uncoupled neuron since $I_{\mathrm{gap},i}(t)=I_{\mathrm{gap},j}(t)=0$ holds at all
times. We use a gap weight $g$ of
$30.0\,\mathrm{nS}$, that represents the typical total coupling of a single neuron
with the remainder of the network: the natural weight of a single
gap junction is much smaller, but each neuron is connected to dozens of other neurons. 
Therefore, the test case exposes how the non-iterative method operates on networks of synchronized neurons coupled by gap junctions.

   \begin{figure}[!tbh]
\begin{centering}
\includegraphics{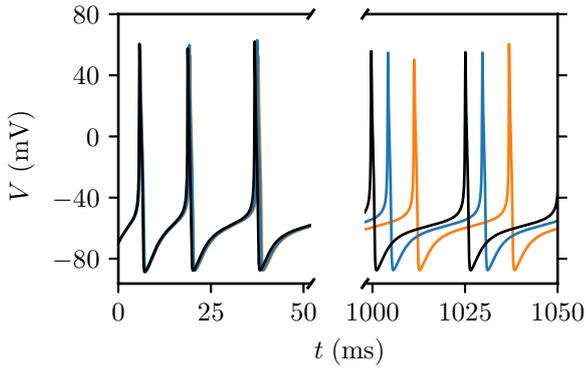}
\par\end{centering}
\vspace*{-0.3cm}
\caption[Artefactual shift without waveform-relaxation techniques]{\textbf{Artefactual shift without waveform-relaxation techniques.} The
black curve shows the reference time course of the membrane potential
of a Hodgkin-Huxley point-neuron model subject to a constant input
current of $200\;\mathrm{pA}$ after simulating $0-50\;\mathrm{ms}$ and $1000-1050\;\mathrm{ms}$. The other curves indicate the time course of the
membrane potential of the same neuron with the same input for the
case that the neuron is coupled by a gap junction to a second model
neuron with exactly the same properties, and the simulation is carried
out with the non-iterative approach using a Runge-Kutta-Fehlberg solver
with an adaptive step-size control to cover the interval of one computational
time step $h=0.1\;\mathrm{ms}$. The orange curve displays the result without spike detection and the blue curve the corresponding results with spike detection.
\label{fig:V_example}}
\end{figure}

\prettyref{fig:V_example} shows that the non-iterative
approach produces a significant shift within only $1\;\mathrm{s}$
of simulated time, when simulated with a commonly used step size, even when spike detection is used. 
This shift results from the calculation of $I_{\mathrm{gap}}$ at intermediate points in an interval of length $h$: 
As the membrane potential of the local neurons evolves over time while the 
membrane potential of the connected neuron is approximated by the extrapolation an artefactual gap current $I_{\mathrm{gap}} \neq 0$ introduces an error to the solution. 
To yield accurate results the time step would have to be exceedingly small requiring significantly more communication.

A suitable waveform-relaxation method can help to overcome these problems. \prettyref{sec:wfr-acc} formulates the specific method and investigates the 
increase of accuracy with respect to the simulation time. Then \prettyref{sec:wfr-per} demonstrates an additional gain of 
the waveform-relaxation method when simulating large-scale networks on supercomputers. All performance benchmarks use an 
implementation in NEST as of version 2.16.0\footnote{The spike detection is not part of this release and was manually added for the investigations in this article.} \citep{nest2160} and run on either a workstation computer or the JURECA 
supercomputer \citep{jureca} at the Juelich Research Centre. 

\subsection{Usage of waveform-relaxation techniques to increase accuracy}\label{sec:wfr-acc}

 For the waveform-relaxation method we consider the ODE-system of each neuron as a subsystem and use the extrapolation of the non-iterative method with spike detection as initial guess $y_{i}^{(0)}(t)$. In order to retain the full potential of parallelization we use a Jacobi waveform-relaxation method, where all neurons can be processed in parallel.
This setup automatically addresses the main accuracy issue of the non-iterative method: Beginning with the second iteration of the waveform-relaxation method the constant or spike-shaped extrapolation of the non-iterative method can be replaced by an interpolation obtained from the results of the previous iteration, which is obviously much closer to the correct solution and improves from iteration to iteration. We choose a cubic Hermite interpolation based on the values of the membrane potential at the grid points and their derivatives that are known through the differential equation. This interpolation allows an straight-forward calculation of the interpolation coefficients that need to be communicated to the other subsystems.

\begin{figure}[!tbh]
\begin{centering}
\includegraphics[width=0.47\textwidth]{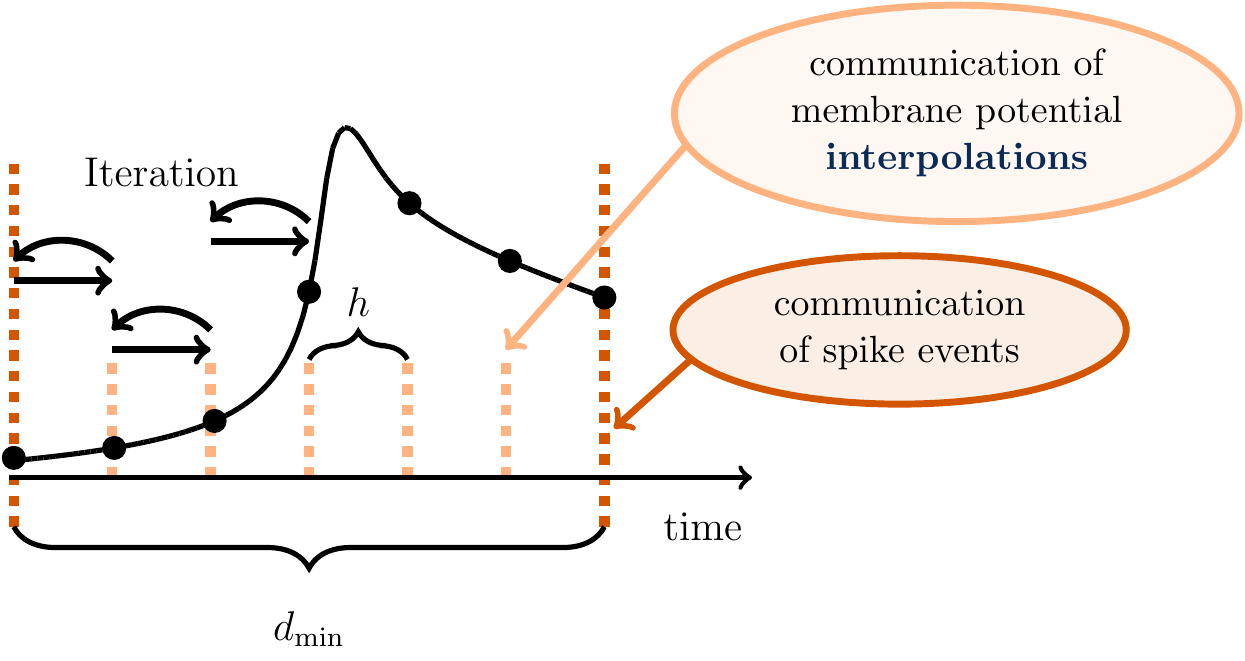}
\par\end{centering}
\caption[Using the waveform-relaxation method to improve the accuracy]{\textbf{Using the waveform-relaxation method to improve the accuracy.} 
The length of the iteration interval of the waveform-relaxation method is chosen as on computation time
step $h$.}\label{fig:com_strategy_acc}
\end{figure}

\prettyref{fig:com_strategy_acc} shows the setup of the waveform-relaxation method from the single-neuron perspective. For now we use the same the same communication points as the non-iterative method (\prettyref{fig:straight-forward}), i.e. we set $\mathcal{T}=h$. In comparison to the non-iterative method the results of each $h$-interval are computed and communicated several times until some convergence criterion is met. 
We stop iterating if

\begin{equation*}
| V_{i}^{(m)}(t_0 + uh)-V_{i}^{(m-1)}(t_0 + uh)|\leq\texttt{wfr\_tol}
\end{equation*}
\noindent holds for of every neuron $i=1,\ldots,v$ and every $u=1,\ldots,\frac{\mathcal{T}}{h}$. Formulated for an arbitrary explicit Runge-Kutta method and an arbitrary length of the iteration interval $\mathcal{T}$ the final methods reads
\begin{equation}
\begin{split}
y_{k+1,i}^{(m)}  \, =  \; & y_{k,i}^{(m)}  + \Delta t_{k+1}^{(m)}\, 
 \sum\limits_{l=1}^{s} b_l  \,f_i  \left(  t^{*(m)}_{k+1,i,l},\right.\, \\[0.15em] & \tilde{y}_{1}^{(m-1)}(t^{*(m)}_{k+1,i,l}),   \ldots,  \tilde{y}_{i-1}^{(m-1)}(t^{*(m)}_{k+1,i,l}), 
\, \\[0.4em] & y_{k+1,i,l}^{*(m)}, \, \tilde{y}_{i+1}^{(m-1)}(t^{*(m)}_{k+1,i,l}),  \ldots,  \left. \tilde{y}_v^{(m-1)}(t^{*(m)}_{k+1,i,l}) \right)
\end{split}
\label{eq:the-ode-method-in-math}
\end{equation}
with
\begin{equation*}
\begin{split}
y^{*(m)}_{k+1,i,q} \, = \; &  y_{k,i}^{(m)} + \Delta t_{k+1}^{(m)}\,
\sum\limits_{l=1}^{q-1} a_{ql}   \,f_i  \left(  t^{*(m)}_{k+1,i,l},\right.\, \\[0.15em] &  \tilde{y}_{1}^{(m-1)}(t^{*(m)}_{k+1,i,l}), \ldots, \tilde{y}_{i-1}^{(m-1)}(t^{*(m)}_{k+1,i,l}),
\\[0.4em] & y_{k+1,i,l}^{*(m)}, \, \tilde{y}_{i+1}^{(m-1)}(t^{*(m)}_{k+1,i,l}), \ldots,\left. \tilde{y}_v^{(m-1)}(t^{*(m)}_{k+1,i,l}) \right)
\end{split}
\end{equation*}
for $q = 1,\ldots,s$ and
\begin{equation*}
\begin{split}
\tilde{y}_i^{(m)}(t_0+(u+\theta) h) \, & = \;  y^{(m)}_{k_u,i}\cdot p_1(\theta) \, + \, y^{(m)}_{k_{u+1},i}\cdot p_2(\theta) \\
& + \; hf(t_0+uh,\, y^{(m)}_{k_u,i})\cdot p_3(\theta) \\
& + \; hf(t_0+(u+1)h,\, y^{(m)}_{k_{u+1},i})\cdot p_4(\theta)
\label{eq:the-limit-method-pre}
\end{split}
\end{equation*}
for $u=0,\ldots,\frac{\mathcal{T}}{h}-1$ and $\theta\in(0,1)$. Here $p_1,p_2,p_3$ and $p_4$ are the Hermite basic polynomials, which read
\begin{eqnarray*}
\; & p_1(\theta) = 1-3\theta^2+2\theta^3, \quad \quad & p_2(\theta)=3\theta^2-2\theta^3,\\
\; & p_3(\theta) = \theta-2\theta^2+\theta^3, \quad \quad \;\; &  p_4(\theta)=-\theta^2+\theta^3
\end{eqnarray*}
and $k_u$ denotes the index that satisfies $t^{(m)}_{k,i}=t_{0} + uh$ for $u=0,\ldots,\frac{\mathcal{T}}{h}$. The existence of those indices is ensured by the used adaptive step-size control.

The method \prettyref{eq:the-ode-method-in-math} is closely related to a general class of waveform-relaxation methods investigated in \citep{Bellen93}.
The only difference that excludes our method from the class of methods discussed in \citep{Bellen93} is the interpolation: The latter work uses interpolation between all grid points, while for us interpolation between  
all grid points $t^{(m)}_{1,i},\ldots,t^{(m)}_{n_{(i,m)},i}$ is neither reasonable, due to the massive amount of communicated data, nor possible, due to the restrictions on the communication points. However, the general convergence of the method to a limit function (for sufficiently small $h$) can be proven analogously. 
The proof (see Theorem 2.2 in \citep{Bellen93}) is based on the contraction principle. The idea is to show that the sequence of mappings
\begin{equation*}
\tilde{y}^{(m)} = \Phi^{(m)}(\tilde{y}^{(m-1)}), \quad m>0, \quad \tilde{y}^{(0)} \text{ given,}
\end{equation*}
is contractive on the space $S$ of all continuous functions $\phi(t)$ of $\left[t_0 , t_0+T \right]$ such that $\phi(t_0)=y_0$. 
In order to show that
\begin{equation*}
\lVert \tilde{a}^{(m)}-\tilde{b}^{(m)} \rVert \;  \overset{!}{<} \; \lVert \tilde{a}^{(m-1)}-\tilde{b}^{(m-1)}  \rVert
\end{equation*}
with $\tilde{a}^{(m-1)},\tilde{b}^{(m-1)} \in S$ holds for some suitable norm, one needs to find an estimate for the difference $\tilde{a}^{(m)}-\tilde{b}^{(m)}$. It is obvious that for our method, this difference can be bounded in the same manner, except that the estimate depends on $h$ instead of the maximal step size $\Delta t^* := \max_{k,i,m} \; \Delta t^{(m)}_{k,i}$ used by the step-size control, as in the original proof (see \citep{Hahne18} for a more detailed comparison). 

We again employ the two-neuron test case to investigate the accuracy of the waveform-relaxation method on the single-neuron level and use the results of an uncoupled neuron with the same properties and inputs 
as reference solution $V_{\mathrm{ref}}$ to determine the quality of the
investigated methods. 

\begin{figure}[!tbh]
\begin{centering}
\includegraphics{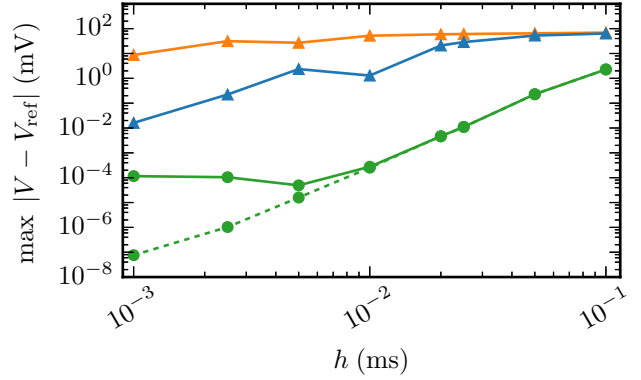}
\par\end{centering}
\caption[Accuracy of the waveform-relaxation method]{\textbf{Accuracy of the waveform-relaxation method.} Error measured as the maximum of $| V - V_{\mathrm{ref}} |$ over all $h$-grid points in $1\;\mathrm{s}$ of simulated time
plotted against the step size $h$. Triangles show results
without waveform-relaxation techniques and without (orange) and with (blue) spike detection. Circles indicate results obtained
by the waveform-relaxation method with cubic interpolation (green).
For the waveform-relaxation method the solid (respectively dashed) curve indicates the choice of the absolute error tolerance of the adaptive step-size mechanism and the \texttt{wfr\_tol} parameter as $10^{-6}$ (respectively $10^{-10}$). \label{fig:results_accuracy}}
\end{figure}

\prettyref{fig:results_accuracy} compares the results of the iterative
method with the results of the non-iterative method without and with spike detection in terms of accuracy
and simulation time. We measure the error $\max \, | V - V_{\mathrm{ref}} |$ of both
methods for different step sizes $h$. For any given step size $h$
the error of the waveform-relaxation method is much smaller than that of
the non-iterative approach, which without spike detection does not even reach a satisfying accuracy
for step size $h=0.001\;\mathrm{ms}$. The error of the waveform-relaxation method for step sizes $h \geq 0.01 \;\mathrm{ms}$ is determined by the interpolation error of the membrane potential, while for smaller step sizes the default error tolerance settings of the neural network simulator prevent further improvement. An adjustment of those settings can, however, restore the original behavior.
In terms of simulation time it is obvious that for any given step size $h$ the non-iterative method
is the fastest implementation, since no additional
iterations are needed to compute the results. In order to evaluate the efficiency of the waveform-relaxation method in comparison to the non-iterative approach we therefore also need to take simulation times in account.

\begin{figure}[!tbh]
\begin{centering}
\includegraphics{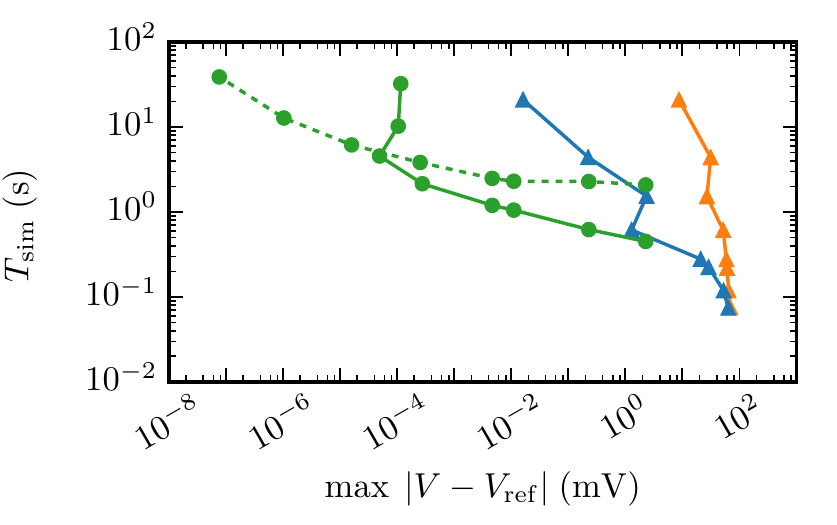}
\par\end{centering}
\caption[Efficiency of the waveform-relaxation method]{\textbf{Efficiency of the waveform-relaxation method.} Simulation time versus error for the
simulations from \prettyref{fig:results_accuracy}.
 \label{fig:results_efficiency}}
\end{figure}

\prettyref{fig:results_efficiency} compares the methods in terms
of efficiency by analyzing the simulation time as a function
of the error. There are two ways of reading this graph: horizontally, one
can find the most accurate method for a given simulation time. Vertically
one can find the fastest method for a desired accuracy. The results
show that the waveform-relaxation method with default error tolerance settings delivers better results in shorter
time than the non-iterative method, even when used with spike detection. As already indicated by \prettyref{fig:results_accuracy} more restrictive error tolerance settings are only beneficial in the unlikely case that the error needs to be reduced below $10^{-4}$ and otherwise only slow down the simulation.

\subsection{Usage of waveform-relaxation techniques to increase performance on supercomputers}\label{sec:wfr-per}

In large-scale simulations with distributed memory communication is expensive, because it is associated with a considerable latency. 
Therefore the total number of communications is an important quantity in large-scale simulations on supercomputers. Although increasing the length of the iteration interval $\mathcal{T}$ slows down
the convergence of waveform-relaxation methods (as indicated by \prettyref{theo:bjorhus-continuous}) and thereby increases the computational load, 
it reduces the total number of communications if $\mathcal{T_{\mathrm{new}}}/\mathcal{T_{\mathrm{old}}}$ is greater than the ratio of the mean number of iterations $\iota_{\mathrm{new}}/ \iota_{\mathrm{old}}$.  
Therefore we investigate if this strategy is able to increase the overall performance in simulations with distributed memory.

\begin{figure}[!tbh]
\begin{centering}
\includegraphics[width=0.47\textwidth]{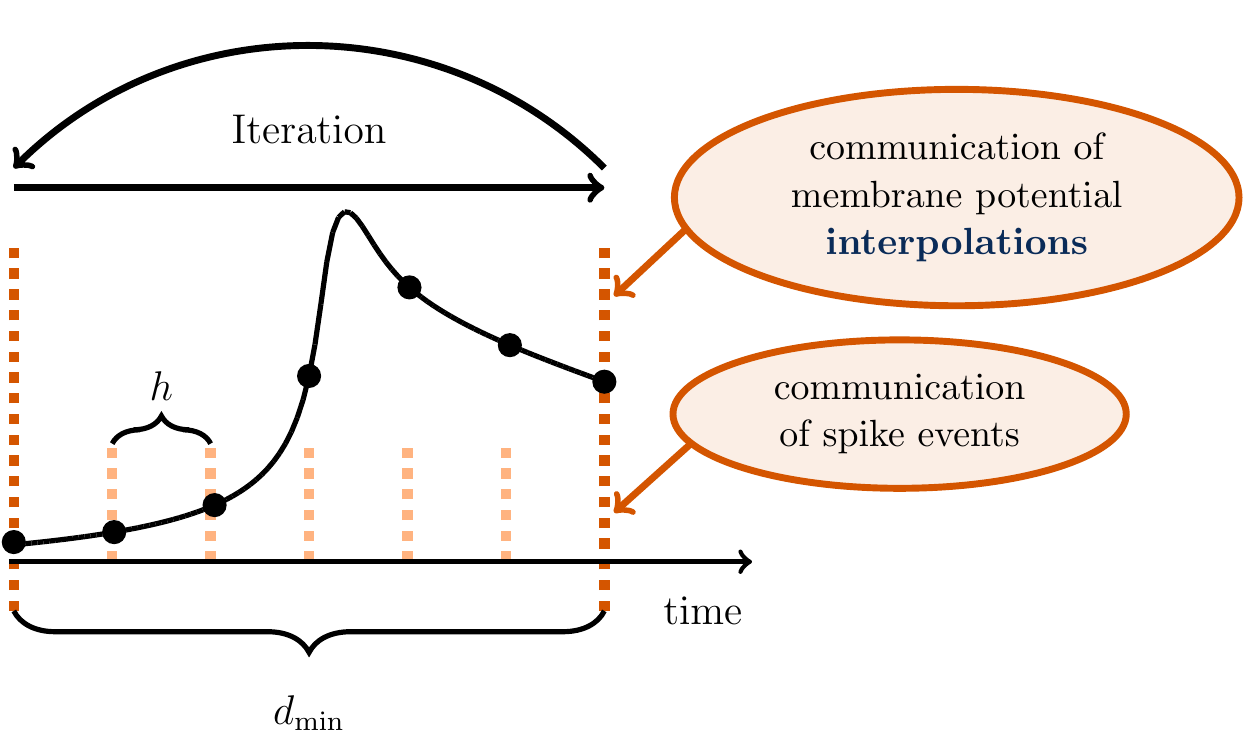}
\par\end{centering}
\caption[Using the waveform-relaxation method to improve the performance on supercomputers]{\textbf{Using the waveform-relaxation method to improve the performance on supercomputers.} 
The length of the iteration interval is chosen as the minimal synaptic delay $d_{\mathrm{min}}$.}  \label{fig:com_strategy_perf}
\end{figure}

\begin{sloppypar}
\prettyref{fig:com_strategy_perf} shows the setup of the waveform-relaxation method with $\mathcal{T}=d_{\mathrm{min}}$ from the single-neuron perspective. Using the minimal network delay as the length of the iteration interval reestablishes the original communication scheme of spiking neural network simulators in the absence of gap junctions (\prettyref{fig:setup-solver}) and therefore constitutes the natural choice for a prolonged iteration interval. The choice between this setup and the one in \prettyref{fig:com_strategy_acc} with $\mathcal{T}=h$ is simply a question of the simulation time, since -- given the convergence of the method -- both strategies deliver nearly identical
results for a given set of parameters.
\end{sloppypar}

\begin{figure}[!tbh]
\begin{centering}
\includegraphics{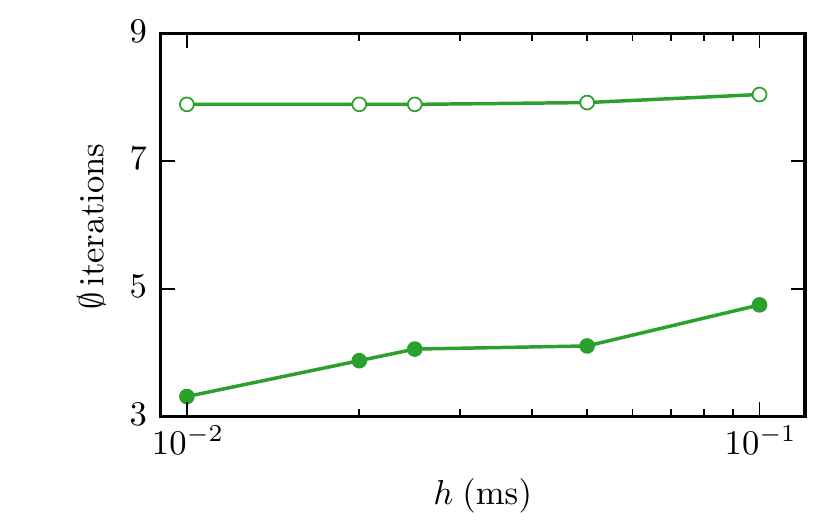}
\par\end{centering}
\caption[Comparison of the mean number of iterations]{\textbf{Comparison of the mean number of iterations.} 
Mean number of iterations measured over $1\,\mathrm{s}$ of simulated time of the two-neuron network for different step sizes $h$ when using the iteration control with $\texttt{wfr\_tol}=10^{-4}$ and a minimal network delay of $d_{\mathrm{min}}=1.0\,\mathrm{ms}$. 
Filled symbols show the results for communication in every step ($\mathcal{T}=h$) while open symbols show the results
for the original NEST communication scheme ($\mathcal{T}=d_{\mathrm{min}}$). As before green circles indicate the waveform-relaxation method with cubic interpolation.
\label{fig:iterations}}
\end{figure}

\begin{sloppypar}
\prettyref{fig:iterations} shows the mean number of iterations for both setups in dependency on the step size $h$ when the minimal network delay is $1.0\,\mathrm{ms}$ in the two-neuron network example. The number of iterations needed is mostly independent of the step size $h$ and differs by about four iterations for the two communication strategies leading to a ratio $\iota_{\mathrm{d_{\mathrm{min}}}}/ \iota_{\mathrm{h}}$ between $1.7$ for $h=0.1$ and $2.5$ for $h=0.01$. At the same time the ratio of communication points $d_{\mathrm{min}}/h$ is much higher, i.e. the overall number of communications reduces significantly. 
\end{sloppypar}

\begin{sloppypar}
We run benchmarks on the JURECA supercomputer in Juelich to investigate the performance of both setups in distributed large-scale simulations.
The benchmarks use
a scaled version of the two-neuron network, where each neuron is coupled to $60$ other neurons by
gap junctions. The number of neurons performing the computation and
the amount of communicated data thus increase with the number of neurons $v$. We keep the conductance
of a neuron accumulated over all gap junctions the same as in our original two-neuron network.
As a consequence, the computations carried out by
each individual neuron are the same, and hence its dynamics is independent
of $v$. Thus, the performance can
be measured in a setting with fixed single neuron dynamics despite
the presence of additional neurons.
\end{sloppypar}

\vspace*{-0.3cm}
\begin{figure}[!tbh]
\begin{centering}
\includegraphics{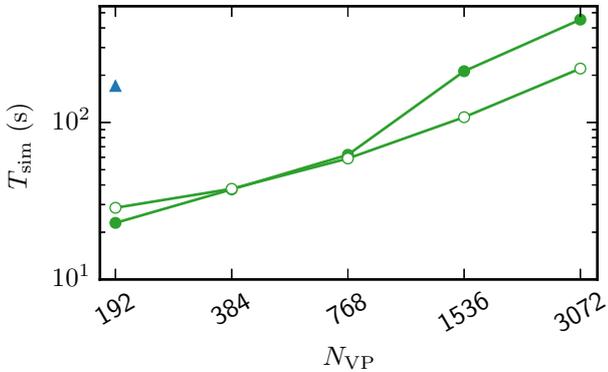}
\par\end{centering}
\caption[Weak scaling on JURECA]{\textbf{Weak scaling on JURECA.}  Simulation of a scaled version of the two-neuron network on JURECA in a weak scaling with $1000$ neurons per compute node and overall $N_{\mathrm{VP}}$ virtual processes ($24$ threads per compute node). Again filled symbols show the results for communication
in every step ($\mathcal{T}=h$) while open symbols show the results
for the original NEST communication scheme ($\mathcal{T}=d_{\mathrm{min}}$). Green circles indicate the waveform-relaxation method with cubic interpolation and step size $h=0.05\,\mathrm{ms}$.
Blue triangles show results without waveform-relaxation techniques and with spike prediction with a step size of $h=0.001\,\mathrm{ms}$.
All simulations use a minimal network delay of $d_{\mathrm{min}}=1.0\,\mathrm{ms}$, a convergence tolerance (\texttt{wfr\_tol}) of $10^{-4}$ and run for $500\,\mathrm{ms}$ of simulated time.
\label{fig:results_jureca_weak}}
\end{figure}

\begin{sloppypar}
\prettyref{fig:results_jureca_weak} shows a weak scaling with $1000$ neurons per compute node. 
Beginning at $64$ compute nodes (each running $24$ threads on the $24$ processors) the waveform-relaxation method 
with $\mathcal{T}=d_{\mathrm{min}}$ reduces the simulation time of the test case significantly, 
while for less compute nodes (and workstation computers, not shown) the original strategy with $\mathcal{T}=h$ results in a better or equal performance.
\end{sloppypar}

\begin{figure}[!tbh]
\begin{centering}
\includegraphics{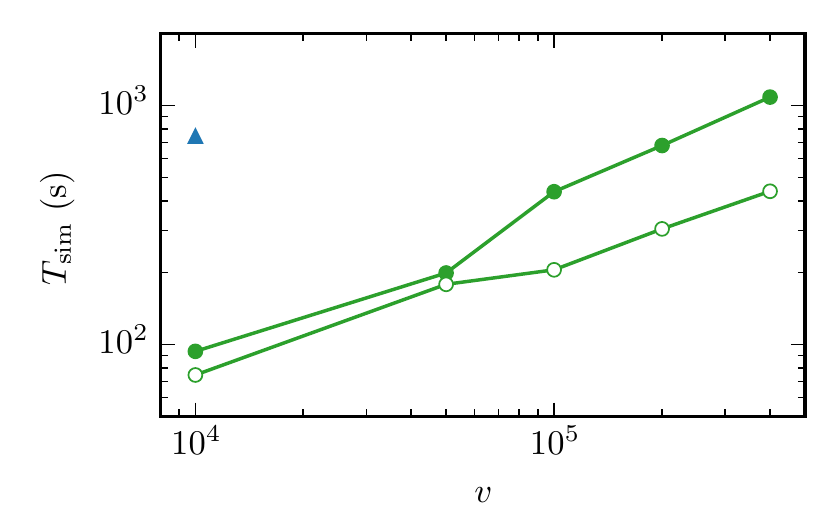}
\par\end{centering}
\caption[Simulations with
different network sizes on JURECA]{\textbf{Simulations with
different network sizes on JURECA.}  Simulation of a scaled version of the two-neuron network on JURECA with
different network sizes $v$ on $3072$ virtual processes ($128$ MPI processes a $24$ threads). All settings and symbols as in \prettyref{fig:results_jureca_weak}.
\label{fig:results_jureca_size}}
\end{figure}

\prettyref{fig:results_jureca_size} demonstrates that the increase in performance for $64$ compute nodes and more can be observed regardless of the network size. 
Both \prettyref{fig:results_jureca_weak} and \prettyref{fig:results_jureca_size} also shows exemplary results for the non-iterative method on JURECA. In this distributed setup the non-iterative method suffers even more from the very small step size $h=0.001$ 
that is necessary to produce an acceptable accuracy.

\section{Conclusion}\label{sec:conclusion}

The inclusion of gap junctions in spiking neural network simulators constitues an interesting use case for waveform-relaxation methods in computational neuroscience. We developed a waveform-relaxation method
that uses a Runge-Kutta method with adaptive step-size control combined with an interpolation on a given coarser grid and is therefore compatible with the general workflow of spiking neural network simulators. The method improves the inclusion of gap junctions in two different ways: First it delivers better results in shorter time than the best possible non-iterative solution, secondly it can be adapted to increase performance for simulations with distributed memory on supercomputers even further. In this work we restricted ourself to a simple two-neuron test case that, however, discloses the behavior of the different methods and is realistic with respect to the conductance of a neuron accumulated over all gap junctions. The shortcomings of the non-iterative method and therefore the need for the proper solution with a waveform-relaxation method can also be shown in more realistic network test cases \citep[Section 3.2]{Hahne15_00022}. The developed waveform-relaxation method still leaves potential for optimization. Replacing the Jacobi iteration scheme and the global convergence condition with more custom-made solutions might further speed-up simulations. However, first each optimization needs to be carefully checked for compability with the rather complex internal structures of spiking neural network simulators.

\addtolength{\textheight}{-2.8cm}   % This command serves to balance the column lengths

%%%%%%%%%%%%%%%%%%%%%%%%%%%%%%%%%%%%%%%%%%%%%%%%%%%%%%%%%%%%%%%%%%%%%%%%%%%%%%%%

%%%%%%%%%%%%%%%%%%%%%%%%%%%%%%%%%%%%%%%%%%%%%%%%%%%%%%%%%%%%%%%%%%%%%%%%%%%%%%%%

%%%%%%%%%%%%%%%%%%%%%%%%%%%%%%%%%%%%%%%%%%%%%%%%%%%%%%%%%%%%%%%%%%%%%%%%%%%%%%%%
%\section*{APPENDIX}


\section*{ACKNOWLEDGMENT}

\begin{sloppypar}The authors gratefully acknowledge the computing time on the supercomputer
JURECA \citep{jureca} at Forschungszentrum J\"{u}lich granted
by the JARA-HPC Vergabegremium (provided on the JARA-HPC partition,
jinb33). This research received funding from the European Union’s Horizon 2020 Framework Programme 
for Research and Innovation under the Framework Partnership Agreement No. 650003 (HBP FPA). 
All spiking network simulations carried out with NEST (http://www.nest-simulator.org).
\end{sloppypar}

\bibliography{the_bib2}

\end{document}